\documentclass[letterpaper, 10 pt, conference]{IEEEconf}  

\IEEEoverridecommandlockouts                              
\usepackage{epsfig}
\usepackage{latexsym}
\usepackage{pdfsync}
\usepackage{amssymb}

\def\Real{{I\!\!R}}

\newcommand{\eq}{\begin{equation}\begin{array}{rllllllllllllllllllllllllllllllll}}
\newcommand{\ee}{\end{array}\end{equation}}
\newcommand{\bmt}{\left[ \begin{array}{ccccccccc}}
\newcommand{\emt}{\end{array}\right]}
\newcommand{\bea}{\begin{eqnarray}}
\newcommand{\eea}{\end{eqnarray}}
\newcommand{\bean}{\begin{eqnarray*}}
\newcommand{\eean}{\end{eqnarray*}}

\title{Adaptive Horizon Model Predictive Control}
\author{Arthur J Krener
\thanks{Research supported in part by AFOSR.}
\thanks{A. J. Krener is with the Department of Applied Mathematics, Naval Postgraduate School, Monterey, CA 93943
        {\tt\small ajkrener@nps.edu}}}%

\begin{document}

\date{}

\maketitle

\begin{abstract}
Adaptive Horizon Model Predictive Control  (AHMPC) is a scheme for varying as needed the horizon length of Model Predictive Control (MPC). Its goal is to achieve stabilization with horizons as small as possible so that MPC can be used on faster or more complicated dynamic processes.  Beside the standard requirements of MPC including a terminal cost that is a control Lyapunov function, AHMPC requires a terminal feedback that turns the control Lyapunov function into a standard Lyapunov function in some domain around the operating point. 
 But this domain need not be known explicitly.  MPC does not compute off-line the optimal cost and 
 the optimal feedback over a large domain instead it computes these quantities on-line when and where they are needed.  AHMPC does not compute off-line the  domain on which the terminal  cost is a control Lyapunov function instead it computes on-line when a state is in this domain. 
   \end{abstract}

\section{Introduction}
\setcounter{equation}{0}
Model Predictive Control (MPC) is a way to optimally steer a discrete time control system to a desired 
operating point.  We briefly describe it following the definitive treatise of Rawlings and Mayne \cite{RM09}.
We  closely follow their notation.

  We are given a controlled, nonlinear dynamics in discrete time 
  \bea \label{dyn}
  x^+&=& f(x,u)
  \eea
  where the state $x\in \Real^{n\times1}$, the control $u\in \Real^{m\times1}$ and $x^+(k)=x(k+1)$.
  This could be the discretization of a controlled, nonlinear dynamics in continuous time.  The goal is 
  to find a feedback law $u(k)=\kappa(x(k))$ that drives the state of the system to some desired 
  operating point.  A pair $(x,u)$ is  an operating point  if $f(x,u)=x$.   We conveniently  assume that the operating point has been translated to be $(x,u)=(0,0)$.    
  
  The controlled dynamics may be subject to constraints such
  as
  \bea \label{sc}
  x&\in& \mathbb{X}\subset \Real^{n\times 1}\\
  u&\in& \mathbb{U}\subset \Real^{m\times 1}
  \eea
  and possibly a constraint involving both the state and control
  \bea \label{scc}
  y=h(x,u)&\in& \mathbb{Y}\subset \Real^{p\times 1}
  \eea
  A control $u$ is said to be feasible at $x\in \mathbb{X}$ if $u\in \mathbb{U}$ and 
   \bean
  f(x,u) & \in &\mathbb{X}\\
   h(x,u) & \in &\mathbb{Y}\\
   \eean
  Of course  the stabilizing  feedback $\kappa(x)$ that we seek needs to be feasible.  For every  $x\in\mathbb{X}$,
  \bean
  \kappa(x)&\in &  \mathbb{U}\\
  f(x,\kappa(x)) & \in &\mathbb{X}\\
   h(x,\kappa(x)) & \in &\mathbb{Y}\\
   \eean

  An ideal way to solve this problem is to choose a Lagrangian $l(x,u)$ that is nonnegative definite 
  in $x,u$ and positive definite in $u$ and then to solve the infinite time optimal control problem of minimizing
  over choice of feasible control sequence ${\bf u}_\infty=(u(0), u(1),\ldots)$ the quantity
  \bea
 V_\infty(x)= \sum_{k=0}^\infty l(x(k),u(k))
  \eea
 subject to the dynamics (\ref{dyn}), the constraints (\ref{sc}, \ref{scc}) and $x(0)=x$.   Let $V^0_\infty(x)$ denote the minimum value 
  and ${\bf u}^0_\infty=(u^0_\infty(0), u^0_\infty(1),\ldots)$ be a minimizing control sequence with corresponding state sequence ${\bf x}^0_\infty=(x^0_\infty(0)=x, x^0_\infty(1),\ldots)$.  Minimizing control and state sequences need not be unique but we shall generally ignore this.
  
  If a pair $V^0_\infty(x)\in \Real, \kappa_\infty(x)\in \Real^{m\times 1}$ of functions satisfy the infinite horizon Dynamic Program Equations (DPE$_\infty$)
  \bea \label{DPEi1}
  V^0_\infty(x)&=&\mbox{min}_u  \left\{ l(x,u)+V^0_\infty(f(x,u))\right\}\\
  \kappa_\infty(x)&=& \mbox{argmin}_u \left\{ l(x,u)+V^0_\infty(f(x,u))\right\}   \label{DPEi2}\\
  V^0_\infty(0)&=&0 \label{DPEi3}
\eea
and the feasibility constraints 
\bea
f(x,\kappa_\infty(x))&\in& \mathbb{X}\\  \label{sc1}
h(x,\kappa_\infty(x))&\in&\mathbb{Y}\label{scc1}
\eea
for $x\in \mathbb{X}$
then it is not hard to show that $V^0_\infty(x)$ is the optimal cost and
$ \kappa_\infty(x)$ is an optimal feedback law,
$
u^0_\infty(k)= \kappa(x^0_\infty(k))
$.
Then under suitable condtions  a Lyapunov argument can be used to show that the feedback
$ \kappa_\infty(x)$ is stabilizing.

The difficulty with this approach is that it is generally impossible to  solve DPE$_\infty$ on a large domain $\mathbb{X}$ if the state dimension $n$ is greater than $2$ or $3$.  So both theorists and practicioners have turned to Model Predictive Control (MPC).  They choose a Lagrangian $l(x,u)$, a horizon length $N$, a terminal domain 
$\mathbb{X}_f\subset\mathbb{X}$ containing $x=0$ and a terminal cost $V_f(x)$ defined and positive definite on  $\mathbb{X}_f$.  Consider the problem of minimizing by choice of feasible ${\bf u}_N=(u_N(0),u_N(1),\ldots,u_N(N-1))$
\bea
V_N(x)&=& \sum_{k=0}^{N-1} l(x(k),u(k)) + V_f(X(N))
\eea
subject to the dynamics (\ref{dyn}), the constraints (\ref{sc}, \ref{scc}), the terminal condition $x(N)\in \mathbb{X}_f$ and the initial condition $x(0)=x$.
Assuming this problem is solvable, let $V^0_N(x)$ denote the optimal cost, 
\bea
V^0_N(x)&=& \min_{{\bf u}_N} V_N(x)
\eea
where the minimum is taken over all feasible ${\bf u}_N$.
Let
${\bf u}_N^0=(u_N^0(0),u_N^0(1),\ldots,u_N^0(N-1))$  and ${\bf x}_N^0=(x_N^0(0)=x,x_N^0(1),\ldots,x_N^0(N))$  denote optimal control and state sequences and  define
\bean
\kappa_N(x)&=& u_N^0(0)
\eean

Let ${\bf X}_N \subset \mathbb{X}$  be defined inductively,
\bean
{\bf X}_0&=& \mathbb{X}_f\\
{\bf X}_{1}&=&\left\{ x\in \mathbb{X}: \exists u\in \mathbb{U}, f(x,u)\in {\bf X}_{0} \land (\ref{scc})
\right\}\\ 
{\bf X}_{N+1}&=&\left\{ x\in \mathbb{X}: \exists u \in \mathbb{U}, f(x,u)\in {\bf X}_{N} \land (\ref{scc})
\right\}
 \eean

The terminal set $\mathbb{X}_f$ is controlled  invariant (aka viable) if for each $\mathbb{X}_f$
there exists a $u\in \mathbb{U}$ such that $f(x,u)\in \mathbb{X}_f$ and the constraints (\ref{scc})
are satisfied.  If this holds then it is not hard to see inductively that the sets 
are nested ${\bf X}_N \subset  {\bf X}_{N+1}$

 If a pair $V^0_N(x), \kappa_N(x)$ defined on ${\bf X}_N$ satisfy the  horizon $N$ Dynamic Program Equations (DPE$_N$)
  \bea \label{DPEN1}
  V^0_N(x)&=&\mbox{min}_u  \left\{ l(x,u)+V^0_N(f(x,u))\right\}\\
  \kappa_N(x)&=& \mbox{argmin}_u \left\{ l(x,u)+V^0_N(f(x,u))\right\}   \label{DPEN2}\\
  V^0_N(x)&=&V_f(x) \mbox{ for }  x \in \mathbb{X}_f\label{DPEN3}
\eea
where the minimum is over all $u\in \mathbb{U}$ that are feasible at $x\in {\bf X}_N$
then it is not hard to show that $V^0_N(x)$ is the optimal cost and
$ \kappa_N(x)$ is an optimal feedback law
$
u^0_N(k)= \kappa(x^0_N(k))
$. If $V_f(x)$ is a control Lyapunov function on $\mathbb{X}_f$ then under suitable conditions a Lyapunov argument can be used to show that the feedback
$ \kappa_N(x)$ is stabilizing on ${\bf X}_N$.   See \cite{RM09} for more details.

As we noted above solving off-line the infinte horizon optimal control problem for all possible states is generally intractable.
The advantage of solving the   horizon $N$ optimal control problem for the current state $x$  is that it possibly can be done on-line as the process evolves.  If the current value of the state is known to be $x\in {\bf X}_N$ then  the finite horizon $N$ optimal control problem is a nonlinear program with finite dimenionsal decision variable ${\bf u}_N\in \Real^{m\times N}$.   If the time step is long enough, if $f,h,l$ are reasonably simple and if $N$ is small enough then this nonlinear program that can be solved in a fraction of one time step for ${\bf u}^0_N$.  Then the first element of this sequence $u_N^0(0)$ is used as the control at the current time.  The system evolves one time step and the process is repeated at the next time.  Conceptually MPC computes an optimal feedback law $\kappa_N(x) =u^0_N(0)$  but only  at  values of $x$ when and where it is needed.

Some authors do away with the terminal cost $V_f(x)$ but there is a theoretical  and  a practical reason to use one.  The theoretical reason is that a control Lyapunov terminal cost facilitates a proof of asymptotic stability via a Lyapunov argument \cite{RM09}.  The practical reason is that one can usually use a shorter horizon $N$ when there is a terminal cost.  A shorter horizon reduces the dimension $mN$ of the decision variables in the nonlinear programs  that need to be solved on-line.  Therefore MPC with a suitable terminal cost can be used for faster and more complicated systems. 

 The ideal terminal cost $V_f(x)$ is $V_\infty(x)$ of the corresponding infinite horizon optimal control  provided that the latter can be accurately computed off-line  on a reasonably large  terminal set $\mathbb{X}_f$. This may be tractable because the terminal set $\mathbb{X}_f$ may be much smaller than $\mathbb{X}$ and only an approximate solution on $\mathbb{X}_f$ may suffice.  For example $V_\infty(x)$ can be locally approximated by the solution of the infinite horizon LQR problem involving the linear part of the dynamics and quadratic part of the Lagrangian at the operating point.

One would expect when the current state $x$ is far from the operating point, a relatively long horizon $N$ is needed to ensure that $x_N^0(N) \in \mathbb{X}_f$  but as the state approaches the operating point   shorter and shorter horizons can be used.  Adaptive Horizon Model Predictive Control (AHMPC) adjusts the horizon of MPC on-line as it is needed.
In the next section we present an ideal version of AHMPC and in the following section we present a practical implementation of AHMPC.  Finally we close with an example.

\section{Ideal Adaptive Model Prediction Control}

We  shall make some standing assumptions.  The first few
are  drawn from Rawlings and Mayne.

 {\it Assumption 1:} (Assumption 2.2 \cite{RM09})\\
 The functions $f(x,u), l(x,u), h(x,u), V_f(x)$ are continuous on some open set
 containing $\mathbb{X} \times \mathbb{U}$,
 $ l(x,u) $ is nonegative definite in $(x,u)$ and positive definite in $u$ 
 on this open set,
 $V_f$ is positive definite on $\mathbb{X}_f$
 and $f(0,0)=0$, $l(0,0)=0$, $V_f(0)=0$.

 {\it Assumption 2:} (Assumption 2.3 \cite{RM09})\\
The sets $\mathbb{X}$ and $\mathbb{X}_f$ are closed, 
  $\mathbb{X}_f\subset \mathbb{X}$, $\mathbb{U}$
  is compact and $\mathbb{X}$, $\mathbb{X}_f$ and $\mathbb{U}$
  contain neighborhoods of their respective origins.

 {\it Assumption 3:} (Assumptions 2.12 and 2.13 of  \cite{RM09})\\
 For all $x\in \mathbb{X}_f$ there exist a feasible $u$ such that
 \bean
 f(x,u)&\in& \mathbb{X}_f\\
 l(x,u)+ V_f(f(x,u)) &\le & V_f(x)
 \eean
This assumption implies that $\mathbb{X}_f$ is controlled invariant
and that $V_f(x)$ is a control Lyapunov function on $\mathbb{X}_f$.

We make some additional assumptions.

{\it Assumption 4:}  For each $x\in \mathbb{X}$ there is 
a nonnnegative integer $N$ and a control
sequence ${\bf u}_N=(u_N(0),\dots ,u_N(N-1))$
such that the corresponding state sequence
${\bf x}_N=(x_N(0)=x,\dots ,x_N(N))$
starting from $x$ satisfies
$x_N(N)\in \mathbb{X}_f$.\\

This assumption allows us to define a function
$N(x)$ on $\mathbb{X}$ as the minimum of all $N$ such that there 
exist such a control
sequence ${\bf u}_N=(u_N(0),\dots ,u_N(N-1))$
and corresponding state sequence
${\bf x}_N=(x_N(0)=x,\dots ,x_N(N))$
starting from $x$  that satisfies
$x_N(N)\in \mathbb{X}_f$.

Then  the nested sets ${\bf X}_N$ defined above are
given by
\bean
{\bf X}_N&=& \left\{ x\in \mathbb{X}: N(x)\le N \right\}
\eean

{\it Assumption 5:}  There exists a nonegative integer $M$
such that $N(x)\le M$ for all $ x\in \mathbb{X}$.  In other words
\bean
{\bf X}_M&=&\mathbb{X}
\eean

These assumptions imply that the usual MPC with horizon
length $M$ is stabilizing on $\mathbb{X}$ by standard arguments
\cite{RM09}.   But it is a waste of time to use horizon
length $M$ when $N(x)$ is substantially smaller $M$. If the current state is $x$, ideal AHMPC  uses
horizon length $N(x)$.  Then as the current state approaches
the terminal set, ideal AHMPC uses shorter and shorter horizons.  When $x$ is 
in the terminal set, ideal AHMPC uses a horizon length of $N=0$.

  In a moment we shall show 
that the function
\bea  \label{VL}
V(x)&=&V^0_{N(x)}(x)
\eea
is a valid Lyapunov function for the closed loop
system which confirms the
stabilizing property of the ideal AHMPC feedback
\bea \label{KL}
\kappa(x)&=& \kappa_{N(x)}(x)
\eea 

But the reason  why this scheme is 
 not practical is that, in general, it is
impossible to compute the function 
$N(x)$.  In the next section we shall
offer a work around but for now 
we study the stabilizing properties
of  ideal AHMPC.  

{\bf Lemma 1}
Assume Assumptions 1-5 hold.  If $N(x)=N $
and if ${\bf u}_N =(u_N(0),\ldots, u_N(N-1))$,
 ${\bf x}_N =(x_N(0)=x,\ldots, x_N(N))$ are a
 control and state trajectory from $x$ such $x(N)\in \mathbb{X}$ 
 then $N(x(k))=N-k$.
  
 {\bf Proof:}   By assumption $x_N^0(N)\in \mathbb{X}_f$  so
 $ x_N^0(N-1)\in {\bf X}_1$, $ x_N^0(N-2)\in {\bf X}_2$, etc.
 So $ N(x_N^0(k))\le N-k$. 
 
 Suppose for some $k$, $x_N^0(k)<N-k$ then $x_N^0(k)\in X_{N-k-1}$ and $x_N^0(k-1)\in X_{N-k-2}$,
 $x_N^0(k-2)\in X_{N-k-3}$, etc. so $x=x_N^0(0)\in X_{N-1}$ which
 contradicts the assumption that  $N(x)=N $.
 $\quad\ \blacksquare$\\

{\bf Lemma 2} (Compare with Lemma 2.14 of \cite{RM09})
Under Assumptions 1-5 then
\bean
V(x^+)\le V(x)-l((x,\kappa(x))
\eean
where $x^+=f(x,\kappa(x))$ and $V$ is defined by (\ref{VL}).\\

{\bf Proof:}  By definition 
\bean
V(x)&=&V^0_{N(x)}(x)\\
&=&\sum_{k=0}^{N(x)-1} l(x^0_{N(x)}(k),u^0_{N(x)}(k))+ V_f(x^0_{N(x)})\\
\eean
where 
${\bf u}^0_{N(x)} =(u^0_{N(x)}(0),\ldots, u_{N(x)}(N(x)-1))$ 
and ${\bf x}_{N(x)} =(x_{N(x)}(0)=x,\ldots, u_{N(x)}(N(x)))$	
are optimizing control and state sequences for the horizon
$N(x)$ optimal control problem so $\kappa(x)=u^0_{N(x)}(0)$
where $\kappa$ is defined by (\ref{KL}).

Let $x^+=f(x,\kappa(x))=x^0_N(x)(1)$,
by Lemma 1, $N(x^+)=N(x)-1$ so 
\bean
V(x^+)= V^0_{N(x)-1}(x^+)
\eean
So 
for any feasible control sequence
$\bar{{\bf u}}=(u(0),\ldots, u(N(x)-2))$ and corresponding state sequence
$\bar{{\bf x}}=(x(0)=x^+,\ldots, x(N(x)-1))$
\bean
V(x^+)\le \sum_{k=0}^{N(x^+)-1} l(x (k),u(k))+ V_f(x{N(x)})
\eean
In particular if we take $\bar{\bf u}=(u^0_{N(x)}(1),\ldots, u_{N(x)}(N(x)-1))$
then $\bar{\bf x}=(x_{N(x)}(1),\ldots, u_{N(x)}(N(x)))$ and
\bean
V(x^+)\le V(x)-l(x,\kappa(x)) \quad\ \blacksquare
\eean
\\

Following Rawlings and Mayne we make the following
assumption.

{\it Assumption 6:} (Assumption 2.16(a) of \cite{RM09})
\\
The stage cost $l $ and the terminal cost $V_f$ satisfy
\bean 
l(x,u) \ge \alpha_1(|x|) && \forall x\in \mathbb{X}, \  \forall u\in \mathbb{U}\\
V_f(x) \le \alpha_2(|x|) && \forall x\in \mathbb{X}_f\\
\eean
where $\alpha_1(\cdot)$ and $\alpha_2(\cdot)$ are ${\it K}_\infty$ functions.

Assumptions 3 and 6 imply that for each $x\in \mathbb{X}_f$  there exists
a feasible $u$ such that 
\bean
V(f(x,u))\le V(x)-\alpha_1(|x|)
\eean

{\bf Proposition 1} (Compare with Proposition 2.17 of \cite{RM09})
\\
(a)  Suppose that Assumptions 1, 2, 3, 4, 5 and 6 are satisfied.
Then there exists  ${\it K}_\infty$ functions $\alpha_1(\cdot)$ and $\alpha_2(\cdot)$ such that
$V(\cdot)$ has the following properties
\bean
V(x)\ge \alpha_1(|x|)&& \forall x\in \mathbb{X}\\
V(x)\le \alpha_2(|x|)&& \forall x\in \mathbb{X}_f\\
V(f(x,\kappa(x)))\le V(x)-\alpha_1(|x|) && \forall x\in \mathbb{X}
\eean

{\bf Proof}  If $N(x)>0$ then the first inequality follows from Assumpion 6(a) and the fact that
$V(x) \ge l(x,\kappa(x))$.  If $N(x)=0$ then the first inequality follows from Assumpion 7(a).
The second inequality follows Assumption 6(a) and the fact that if $x\in \mathbb{X}_f$ then $N(x)=0$ and
so $V(x)=V_f(x). \quad\ \blacksquare$

  If the second property held for all $x\in \mathbb{X}$
  \bean
V(x)\le \alpha_2(|x|)&& \forall x\in \mathbb{X}
\eean
then  $V(x)$ would be a valid Lyapunov on $\mathbb{X}$.  The following is a paraphrase
of  a proposition  of Rawlings and Mayne.

{\bf Proposition 2} (Proposition 2.1 of \cite{RM09})
\\
Suppose that Assumptions 1, 2, 3 hold, that $\mathbb{X}_f$ contains an open neighborhood of the origin
and that $\mathbb{X}$ is compact.  If there exists a ${\it K}_\infty$ function $\alpha(\cdot)$ such that 
$V^0_N(x)\le \alpha(|x|)$ for $x\in \mathbb{X}_f$ then there exists another ${\it K}_\infty$ function $\beta_N(\cdot)$ such that $V^0_N(x)\le \beta_N(|x|)$  for $x\in {\bf X}_N$.

This allows us to prove the following proposition.

{\bf Proposition 3}\\
Suppose that Assumptions 1, 2, 3, 4, 5 and 6 hold, that $\mathbb{X}_f$ contains an open neighborhood of the origin
and that $\mathbb{X}$ is compact.  If there exists a ${\it K}_\infty$ function $\alpha(\cdot)$ such that 
$V^0_N(x)\le \alpha(|x|)$ for $x\in \mathbb{X}_f$ then there exists another ${\it K}_\infty$ function $\beta(\cdot)$ such that $V^0_N(x)\le \beta(|x|)$  for $x\in \mathbb{X}$.

{\bf Proof:} By Assumpions 4 and 5, $\mathbb{X}={\bf X}_M$.  Let $\beta_0(\cdot)=\alpha(\cdot)$ and define
\bean 
\beta(s)&=& \max \left\{ \beta_N(s): N=0,1,\ldots,M\right\}
\eean
The maximum of a finite family of ${\it K}_\infty$ functions is also a ${\it K}_\infty$ function.   Clearly if 
$x\in \mathbb{X}={\bf X}_M$
\bean
V(x)=V^0_{N(x)}(x) \le \beta_{N(x)}(|x|)\le \beta(|x|) \quad\  \blacksquare
\eean

{\bf Proposition 4} Suppose that Assumptions 1, 2, 3, 4, 5 and 6 hold, that $\mathbb{X}_f$ contains an open neighborhood of the origin and that $\mathbb{X}$ is compact.  Then $V(x)$ is a valid Lyapunov function which confirms the asymptotic stability of the closed loop dynamics
\bean
x^+&=& f(x,\kappa(x)
\eean
 on $\mathbb{X}$

So ideal AHMPC solves our stabilization problem in theory.  But it generally can't be implemented because we can't compute the key ingredient, the function $N(x)$ or its domain of definition.  \\

There is a slightly less ideal version of AHMPC.
  Suppose we have a function $N(x)$ with the following properties. \\
a)   For each $x\in \mathbb{X}$ there is a feasible control sequence ${\bf u}_{N(x)}=(u_{N(x)}(0),\ldots,u_{N(x)}(N(x)-1))$  and corresponding state sequence ${\bf x}_{N(x)}=(x_{N(x)}(0)=x,\ldots,x_{N(x)}(N))$ such that $x_{N(x)}(N(x)))\in \mathbb{X}_f$.
\\
b)  There exist an $M$ such that $N(x)\le M$ for all $x\in \mathbb{X}$.\\
c)  If ${\bf u}^0_{N(x)}=(u^0_{N(x)}(0),\ldots,u^0_{N(x)}(N(x)-1))$  and  ${\bf x}^0_{N(x)}=(x^0_{N(x)}(0)=x,\ldots,x^0_{N(x)}(N))$  are optimal control and state sequences for the horizon $N(x)$ optimal control problem starting at $x$ then $0\le N(x_{N(x)}(k))-N(x_{N(x)}(k+1))\le 1$.  In other words
along optimal trajectories $N(\cdot)$ either stays the same or decreases by $1$ at each time step.

Then the above results hold, $V(x)$ as defined by (\ref{VL})  is a valid Lyapunov function for closed loop system using feedback $\kappa(x)$ defined by (\ref{KL}).    The only additional thing that needs to shown is that if $N(x )=N(x^+)$ then $ V(x)-V(x^+)\ge \alpha_1(|x|)$.  But this follows from standard MPC arguments, see Lemma 2.14 of \cite{RM09}.

\section{Adaptive Horizon Model Predictive Control}
Here is a variation on the above that is practical which we call Adaptive Horizon Model Predictive Control (AHMPC).  
We assume that we have the following.
\begin{enumerate}
\item Sets $\mathbb{X},\  \mathbb{X}_f, \ \mathbb{U}$ satisfying Assumption 2.  We do not require that $ \mathbb{X}_f$ be known explicitly.
\item A discrete time controlled dynamics $f(x,u)$,  a Lagrangian  $l(x,u)$, a constraint pair $(h(x,u),\ \mathbb{Y})$   and a terminal cost $V_f(x)$ satisfying Assumption 1. 
\item A terminal feedback $u=\kappa_f(x)$ and a class ${\it K}_\infty$ function $\alpha(\cdot)$ defined
for all $ x\in  \mathbb{X}_f$ and satisfying
\bean
V_f(x)&\ge& \alpha(|x|)\\
f(x,\kappa_f(x))&\in&  \mathbb{X}_f\\
V_f(x)-V_f(f(x,\kappa_f(x)))&\ge &\alpha(|x|)\\
h(x,\kappa_f(x))&\in& \mathbb{Y}
\eean
\end{enumerate}
We don't need to know the terminal set $\mathbb{X}_f$ on which these conditions are satisfied, all we need to there is such a terminal set and that it contains a neighborhood of $x=0$.

One way of obtaining such a terminal pair $V_f(x), \ \kappa_f(x)$ is to approximately solve the infinite horizon dynamic program equations (DPE$_\infty$) on some neighborhood of the origin.  For example if the linear part of the dynamics and the quadratic part of the Lagrangian constitute a nice LQR problem then then one can let 
$V_f(x)$ be the quadratic optimal cost and $\kappa_f(x)$ be the linear optimal feedback of the LQR.  Alternatively one can take higher degree Al'brekht approximations to $V_\infty(x), \ \kappa_\infty(x)$
\cite{Al61}.  Of course the problem with such terminal pairs $V_f(x), \ \kappa_f(x)$ is that generally there is no way to estimate the terminal set $\mathbb{X}_f$ on which (1), (2) and (3) are satisfied.  It is reasonable to expect that they are satisfied on some terminal set but the extent of the terminal set is very difficult to estimate.

AHMPC mitigates this difficulty.  MPC does not try to compute the optimal cost and optimal feedback everywhere, instead it computes them just when and where they are needed.  AHMPC does not try to compute the extent of $V_f$, it just tries to determine if the  end state $x_N^0(N)$ of the currently computed optimal trajectory  is in a terminal  set $\mathbb{X}_f$ where  (1), (2) and (3) are satisfied.  
 
 Suppose the current state is $x$ and we have solved the horizon $N$ optimal control problem for
 ${\bf u}^0_N=(u^0_N(0),\ldots, u^0_N(N-1))$,  ${\bf x}^0_N=(x^0_N(0)=x,\ldots, x^0_N(N))$. AHMPC does not explictly impose 
 the terminal constraint $x^0(N)\in \mathbb{X}_f$ because $\mathbb{X}_f$ is not explicitly known but it does require that the terminal cost $V_f$ is defined at $x^0(N)$.
 
 The terminal 
 feedback $u=\kappa_f(x)$ is used to extend the state trajectory $L$ additional steps 
  \bean
 x^0_N( k+1))&=& f(x^0_N(k),\kappa_f(x^0_N(k))
 \eean 
for $k=N,\ldots, N+L-1$.   This assumes that the terminal feedback is defined $u=\kappa_f(x)$ is defined on $x^0_N(k)$ for $k=N,\ldots, N+L-1$.  If the terminal  feedback is not defined at any of these points then we presume that $ x^0_N( N) $ is not in $\mathbb{X}_f$ so we increase $N$ by $1$ and we solve the optimal control problem over the new horizon. 

If the feedback is defined on the extended trajectory then one checks that the Lyapunov conditions hold for the extended part of the state sequence,
\bea \label{L1}
V_f (x^0_N(k)&\ge & \alpha(|x^0_N(k)|)\\
V_f (x^0_N(k)-V_f (x^0_N(k+1)&\ge & \alpha(|x^0_N(k)|)
\label{L2}
\eea
for $k=N,\ldots, N+L-1$.   Again if the terminal cost $V_f$ is not defined at any of these points then we presume that $ x^0_N( N) $ is not in $\mathbb{X}_f$ so we increase $N$ by $1$ and we solve the optimal control problem over the new horizon.

If (\ref{L1}, \ref{L2}) hold for all for $k=N,\ldots, N+L-1$.  then we presume that $x^0_N(N)\in \mathbb{X}_f$
and we use the control $u_N^0(0)$ to move one time step forward to  $x^+=f(x, u^0_N(0))$.  At this next state $x^+$ we 
solve the horizon $N-1$ optimal control problem and check that the extension of the new optimal trajectory satisfies (\ref{L1}, \ref{L2}).

If (\ref{L1}, \ref{L2}) do not hold for all for $k=N,\ldots, N+L-1$.  then we presume that $x^0_N(N)\notin \mathbb{X}_f$.
If time permits we solve the horizon $N+1$ optimal control problem at the current state $x$ and then check the Lyapunov  conditions (\ref{L1}, \ref{L2}) again.  We keep increasing the horizon by $1$ until these conditions are satisfied.
If we run out of time before (\ref{L1}, \ref{L2}) are satisfied then we use the last computed $u_N^0(0)$ and move one time step forward to  $x^+=f(x, u^0_N(0))$.  At $x^+$  we 
solve the horizon $N+1$ optimal control problem.

The number $L$ of additional time steps is a design parameter.  Two obvious choices are to take a fixed $L$ which is a fraction of $M$ or to take a varying $L$ which is a fraction of the current $N$.

\section{Example}
The example that we apply AHMPC to is stabilizing a double pendulum to the upright position using torques at each of the pivots.  The states are  $x_1$, the angle of the first leg  measured in radians counter-clockwise from straight up, $x_2$, the angle of the second leg measured in radians counter-clockwise from straight up, $x_3=\dot{x}_1$ and $x_4=\dot{x}_2$.  The controls are $u_1$, the torque applied at the base of the first leg, and $u_2$, the torque applied at the joint between the  legs.  The length of the first leg is $1$ m. and the length of the second leg is $2$ m.  The legs are assumed to  be massless but there is a mass of $2$ kg. at the joint between the legs and a mass of $1$ kg. at the tip of the second leg.   The continuous time controlled dynamics is  discretized using Euler's method with time step $0.1$ s. assuming the control is constant throughout the 
the time step.

The continuous time Lagrangian is chosen to be $l_c(x,u)=(|x|^2+|u|^2)/2$ and its Euler discretization, $l(x,u)=(|x|^2+|u|^2)/20$, is used.  We choose the initial state to be $x=(\pi/2,-\pi/2,0,0)'$  and the initial horizon length to be
$N=5$.  We simulated practical AHMPC with $V_f(x), \ \kappa_f(x)$ the solution of the LQR problem using the linear part of the dynamics at the origin and the quadratic Lagrangian, $\alpha(|x|)=0.1|x|^2$ and fixed $L=5$.  We did not move one time step forward if (\ref{L1}, \ref{L2}) did not hold over the extended state trajectory but instead increased $N$ by one and recomputed.   The AHMPC trajectories of the two angles, $x_1$ in blue and $x_2$ in red,  are shown in Figure 1.
\begin{figure}
\includegraphics[width=3in]{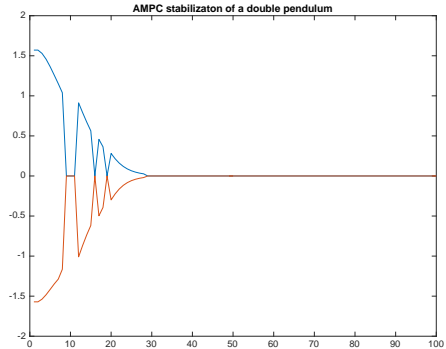} 
\vspace{-1.1in}
\caption{Angles  Converging to the Vertical}    
 \end{figure}
The adaptively changing horizon length is shown in Figure 2.  This graph includes cases where the horizon was increased by one but the state of the pendulum was not advanced.   Notice that the horizon goes down and up several times before settling at $N=0$.
\begin{figure}
\includegraphics[width=3in]{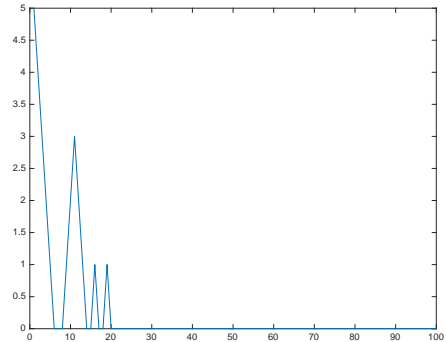}     
\vspace{-1.1in}
\caption{Adaptively Changing Horizon}
 \end{figure}

 \section{Conclusion}  
 Adaptive Horizon Model Predictive Control is a scheme for varying the horizon length  in Model Predictive Control as the stabilization process evolves.
We have presented an ideal  version  of AHMPC and shown that it guarantees stabilization.  AHMPC is a practical version that proceeds without knowing the minimum horizon length function $N(x)$  and without knowing the domain of Lyapunov stability of the terminal cost $V_f(x)$ and terminal feedback $\kappa_f(x)$.

We have only proven the convergence of AHMPC under ideal conditions but
the convergence of standard MPC is also proven under similar ideal conditions, e.g., exact model, exact knowledge of the current state, exact solution of the finite horizon optimal control problems, etc.

The principal advantge of AHMPC over standard MPC is that   the AHMPC horizon length decreases as the process is stabilized thereby lessening the on-line computational burden.  Hence AHMPC may be able to stabilize systems with faster or more complicated dynamics.

The author would like to acknowledge helpful communications with Sergio Lucia,  Philipp Rumschinski  
and Rolf Findeisen.

\end{document}